\documentclass{article}
\usepackage{amssymb}
\usepackage{amsthm}
\usepackage{amsmath}
\usepackage{hyperref}
\usepackage{eucal}
\usepackage[dvips]{graphicx}
\usepackage{fancyhdr}

\makeatother  
\begin{document}
\setlength\textheight{6.9in} \setlength\textwidth{5in}
\setlength\baselineskip{13pt}


\title{A Note on the Smoluchowski-Kramers Approximation for the Langevin Equation with Reflection}

\author{Konstantinos Spiliopoulos\\
Department of Mathematics, University of Maryland\\
College Park, 20742, Maryland, USA\\
kspiliop@math.umd.edu}

\vspace{2cm}
\date{}

\maketitle

\begin{abstract}
According to the Smoluchowski-Kramers approximation, the solution
of the equation
$\mu\ddot{q}^{\mu}_{t}=b(q^{\mu}_{t})-\dot{q}^{\mu}_{t}+\Sigma(q^{\mu}_{t})\dot{W}_{t},
q^{\mu}_{0}=q, \dot{q}^{\mu}_{0}=p$ converges to the solution of
the equation $\dot{q}_{t}=b(q_{t})+\Sigma(q_{t})\dot{W}_{t},
q_{0}=q$ as $\mu\rightarrow 0$. We consider here a similar result
for the Langevin process with elastic reflection on the boundary.
\bigskip

\textit{Keywords}: Smoluchowski-Kramers approximation, reflection,
Langevin equation, Skorohod reflection problem.
\end{abstract}


\vspace{2cm}

\framebox{\parbox[c]{11cm}{This is an electronic reprint of the
original article published by the
\href{http://www.worldscinet.com/sd/sd.shtml}{World Scientific
Publishing Company} in Stochastics and Dynamics, Vol. 7, No. 2, June
2007 , 141-153. This reprint differs from the original in pagination
and typographic detail.}}


\newpage

\section{Introduction}

The well-known Smoluchowski-Kramers approximation
(\cite{K5},\cite{K4}) implies that the solution of the stochastic
differential equation (S.D.E.)
\begin{eqnarray}
\mu\ddot{q}^{\mu}_{t}&=&b(q^{\mu}_{t})-\dot{q}^{\mu}_{t}+\Sigma(q^{\mu}_{t})\dot{W}_{t}\label{sk1}\\
q^{\mu}_{0}&=&q\in\mathbb{R}^{r} \nonumber\\
\dot{q}^{\mu}_{0}&=&p\in\mathbb{R}^{r} \nonumber
\end{eqnarray}
converges in probability as $\mu\rightarrow 0$ to the solution of
the following S.D.E.:
\begin{eqnarray}
\dot{q}_{t} &=& b(q_{t})+\Sigma(q_{t})\dot{W}_{t}\label{sk2}\\
q_{0}&=&q\in\mathbb{R}^{r}, \nonumber
\end{eqnarray}
where $b=(b_{1},...,b_{r})^{'}$ (the transpose of
$(b_{1},...,b_{r})$) with $b_{j}:\mathbb{R}^{r}\rightarrow
\mathbb{R},j=1,..,r$, $\Sigma=[\sigma_{ij}]_{i,j}^{r}$ with
$\sigma_{ij}:\mathbb{R}^{r}\rightarrow \mathbb{R},i,j=1,..,r$ have
bounded first derivatives and
$W_{t}=(W_{t}^{1},...,W_{t}^{r})^{'}$ is the standard
r-dimensional Wiener process. In other words, one can prove that
for any $\delta,T>0$ and $q,p\in\mathbb{R}^{r}$ (see, for example,
Lemma 1 in \cite{K2}),
\begin{eqnarray}
\lim_{\mu\downarrow 0}P(\max_{0\leq t \leq
T}|q^{\mu}_{t}-q_{t}|>\delta)=0. \label{convergence}
\end{eqnarray}

Equation (\ref{sk1}) describes the motion of a particle of mass
$\mu$ in a force field $b(q)+\Sigma(q)\dot{W}_{t}$ with a friction
proportional to velocity. The Smoluchowski-Kramers approximation
justifies the use of equation (\ref{sk2}) to describe the motion
of a small particle.

\setlength\textheight{7.7in} \setlength\textwidth{5in}

It is easy to see now that (\ref{sk1}) can be equivalently written
as:
\begin{eqnarray}
\dot{q}^{\mu}_{t}&=&p^{\mu}_{t}\nonumber\\
\mu\dot{p}^{\mu}_{t}&=&b(q^{\mu}_{t})-p^{\mu}_{t}+\Sigma(q^{\mu}_{t})\dot{W}_{t}\label{sk1system1}\\
q^{\mu}_{0}&=&q\in\mathbb{R}^{r}, \hspace{0.1cm}
\dot{q}^{\mu}_{0}=p\in\mathbb{R}^{r}.\nonumber
\end{eqnarray}

Let us define $\mathbb{R_{+}}=\{q^{1} \in \mathbb{R}: q^{1}\geq
0\}$ and let the configuration space be $D=\mathbb{R_{+}}\times
\mathbb{R}^{r-1}$. In this paper we examine the behavior of the
process with elastic reflection on the boundary $\partial D\times
\mathbb{R}^{r}=(\partial \mathbb{R_{+}}\times
\mathbb{R}^{r-1})\times \mathbb{R}^{r}$ of the phase space
$D\times \mathbb{R}^{r}$ that is governed by (\ref{sk1system1}),
i.e. of the Langevin process with reflection, as $\mu\rightarrow
0$ when $\Sigma$ is the unit matrix. We will show that the first
component (the q component) of the Langevin process with
reflection at $q^{1}=0$, that is governed by equation
(\ref{sk1system1}), converges in distribution to the diffusion
process with reflection on $\partial D$ that is governed by
(\ref{sk2}). The method is based on properties of the Skorohod
reflection problem and in techniques developed in \cite{C1} and in
\cite{C2}. In section 2 we define the Langevin process with
reflection for general diffusion matrx $\Sigma$ with inputs that
have bounded first derivatives, in section 3 we describe the
Skorohod reflection problem and in section 4 we consider the limit
$\mu\rightarrow 0$ when the diffusion matrix is the unit matrix.
We note here that the limit when $\mu\rightarrow 0$ for a general
diffusion matrix as above can be examined similarly.

\section{Langevin process with reflection and preliminary results}
We begin with the construction of the Langevin process
$(q^{\mu}_{t};p^{\mu}_{t})$ in $D\times \mathbb{R}^{r}$ with
elastic reflection on the boundary.  Let $b=(b_{1},...,b_{r})^{'}$
with $b_{j}:D\rightarrow \mathbb{R},j=1,..,r$ and
$\Sigma=[\sigma_{ij}]$ with $\sigma_{ij}:D\rightarrow
\mathbb{R},i,j=1,..,r$ have bounded first derivatives and $\Sigma$
be non-degenerate. Let $(q,p) \in D\times \mathbb{R}^{r}$ be the
initial point (we assume that $(q^{1})^{2}+(p^{1})^{2}\neq 0$).
Then $(q^{\mu}_{t};p^{\mu}_{t})$ is the right-continuous Markov
process in $D\times \mathbb{R}^{r}$ defined as follows. Consider
the following system of S.D.E.'s:
\begin{eqnarray}
\dot{q}^{i,\mu}_{t}&=& p^{i,\mu}_{t}\nonumber \\
\mu\dot{p}^{i,\mu}_{t}&=&-p^{i,\mu}_{t}+b_{i}(q^{\mu}_{t})+\sum_{j=1}^{r}\sigma_{ij}(q^{\mu}_{t})\dot{W}_{t}^{j}\label{sk1system}\\
q^{i,\mu}_{0}&=&q^{i}, \hspace{0.1cm} p^{i,\mu}_{0}=
p^{i},\hspace{0.1cm} i=1,...,r. \nonumber
\end{eqnarray}

We define $(q^{\mu}_{t};p^{\mu}_{t})$ to be the solution to
(\ref{sk1system}) for $t \in [0,\tau_{1}^{\mu})$, where
$\tau_{1}^{\mu}=\inf\{t>0:q_{t}^{1,\mu}=0\}$. Then define
$(q^{\mu}_{t};p^{\mu}_{t})$ for $t \in
[\tau_{1}^{\mu},\tau_{2}^{\mu})$, where
$\tau_{2}^{\mu}=\inf\{t>\tau_{1}^{\mu}:q_{t}^{\mu}=0\}$, to be the
solution of (\ref{sk1system}) with initial conditions
$$(q_{\tau_{1}^{\mu}}^{\mu};p_{\tau_{1}^{\mu}}^{\mu})=(0,\lim_{t
\uparrow \tau_{1}^{\mu}}q_{t}^{2,\mu},...,\lim_{t \uparrow
\tau_{1}^{\mu}}q_{t}^{r,\mu};-\lim_{t \uparrow
\tau_{1}^{\mu}}p_{t}^{1,\mu},\lim_{t \uparrow
\tau_{1}^{\mu}}p_{t}^{2,\mu},...,\lim_{t \uparrow
\tau_{1}^{\mu}}p_{t}^{r,\mu}).$$ If
$0<\tau_{1}^{\mu}<\tau_{2}^{\mu}<...<\tau_{k}^{\mu}$ and
$(q^{\mu}_{t};p^{\mu}_{t})$ for $t \in [0,\tau_{k}^{\mu})$ are
already defined, then define $(q^{\mu}_{t};p^{\mu}_{t})$ for $t
\in [\tau_{k}^{\mu},\tau_{k+1}^{\mu})$ as solution of
(\ref{sk1system}) with initial conditions
$$(q_{\tau_{k}^{\mu}}^{\mu};p_{\tau_{k}^{\mu}}^{\mu})=(0,\lim_{t \uparrow
\tau_{k}^{\mu}}q_{t}^{2,\mu},...,\lim_{t \uparrow
\tau_{k}^{\mu}}q_{t}^{r,\mu};-\lim_{t \uparrow
\tau_{k}^{\mu}}p_{t}^{1,\mu},\lim_{t \uparrow
\tau_{k}^{\mu}}p_{t}^{2,\mu},...,\lim_{t \uparrow
\tau_{k}^{\mu}}p_{t}^{r,\mu})$$ (see Figure 1 for an
illustration).

This construction defines the process $(q^{\mu}_{t};p^{\mu}_{t})$
in $D\times \mathbb{R}^{r}$ for all $t\geq 0$. This follows from
Theorem 2.4, which states that the process that we constructed
above does not have infinitely many jumps in any finite time
interval $[0,T]$. Therefore we have the following definition:

\textbf{Definition 2.1.} \textit{We call the above recursively
constructed process, the Langevin process with elastic reflection
on the boundary $\partial D\times \mathbb{R}^{r}$. This process
has jumps on $\partial D\times \mathbb{R}^{r}$ and is continuous
inside $D\times \mathbb{R}^{r}$.}

We will refer to the Langevin process with reflection as
l.p.r.$(q_{t}^{\mu};p_{t}^{\mu})$. Moreover we will denote by
$(q^{\mu,q}_{t};p^{\mu,p}_{t})$ the trajectories of
$(q_{t}^{\mu};p_{t}^{\mu})$ with initial position $(q,p)$. For
easy of notation we also define $-x=(-x^{1},x^{2},\dots,x^{r})$
and $|x|=(|x^{1}|,x^{2},\dots,x^{r})$ for $x \in \mathbb{R}^{r}$.

Below we see an illustration of the construction above in the
$(q^{1}-p^{1})$ phase space.

\begin{figure}[ht]
\begin{center}
\includegraphics[scale=1, width=7 cm, height=6 cm]{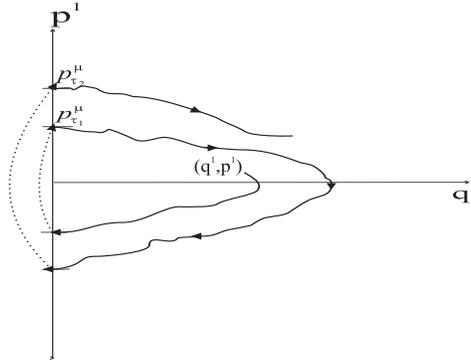}
\caption{Illustration of the Langevin process with reflection in
the $(q^{1}-p^{1})$ phase space}
\end{center}
\end{figure}
Let us give now another construction of the Langevin process with
reflection. Consider the following S.D.E. in $\mathbb{R}^{2r}$:
\begin{eqnarray}
\dot{q}^{1,\mu}_{t}&=& p^{1,\mu}_{t}\nonumber \\
\mu\dot{q}^{1,\mu}_{t}&=&-p^{1,\mu}_{t}+ \textrm{sgn}(q^{1,\mu}_{t})b_{1}(|q^{\mu}_{t}|)+\sum_{j=1}^{r}\textrm{sgn}(q^{1,\mu}_{t})\sigma_{1j}(|q^{\mu}_{t}|)\dot{W}_{t}^{j}\nonumber\\
q^{1,\mu}_{0}&=&q^{1}, p^{1,\mu}_{0}= p^{1}, \nonumber \\ \nonumber \\
\dot{q}^{i,\mu}_{t}&=& p^{i,\mu}_{t}\label{sk1_sgn_withb} \\
\mu\dot{p}^{i,\mu}_{t}&=&-p^{i,\mu}_{t}+ b_{i}(|q^{\mu}_{t}|)+\sum_{j=1}^{r}\sigma_{ij}(|q^{\mu}_{t}|)\dot{W}_{t}^{j}\nonumber\\
q^{i,\mu}_{0}&=&q^{i}, p^{i,\mu}_{0}= p^{i},i=2,...,r, \nonumber
\end{eqnarray}
where $\textrm{sgn}(x)$ takes two values, 1 if $x\geq 0$ and -1 if
$x<0$.

\vspace*{0.2cm} \textbf{Lemma 2.2.} \textit{Equation
(\ref{sk1_sgn_withb}) has a weak solution which is unique in the
sense of probability law.}

\textbf{Proof.} The existence follows from the Girsanov's Theorem
on the absolute continuous change of measures in the space of
trajectories (b and $\Sigma$ are assumed bounded) and the fact
that (\ref{sk1_sgn_withb}) with $b=0$ has a weak solution. The
uniqueness follows from Proposition 5.3.10 of \cite{K3}.

\begin{flushright}
$\Box$
\end{flushright}

Using the processes $(q^{\mu,q}_{t};p^{\mu,p}_{t})$ and
$(q^{\mu,-q}_{t};p^{\mu,-p}_{t})$ we can give another construction
of the Langevin process with reflection, as follows. Assume that
$q^{1}>0$ and $p^{1}>0$, The graphs of $p^{1,\mu,p^{1}}_{t}$ and
$p^{1,\mu,-p^{1}}_{t}$ will be exactly symmetric with respect to
zero. The same will be true also for the graphs of
$q^{1,\mu,q^{1}}_{t}$ and of $q^{1,\mu,-q^{1}}_{t}$. Let
$\tau_{0}^{\mu}=0,\tau_{k}^{\mu}=\inf\{t>\tau_{k-1}^{\mu}:q_{t}^{1,\mu,q^{1}}=0\}$
and $(\widehat{q}^{\mu}_{t};\widehat{p}^{\mu}_{t})$ be a
stochastic process, which is defined as follows:
\begin{eqnarray}
(\widehat{q}^{\mu}_{t};\widehat{p}^{\mu}_{t})&=&(q^{\mu,q}_{t};p^{\mu,p}_{t})
\textrm{ for } \tau_{2k}^{\mu}\leq t\leq
\tau_{2k+1}^{\mu,-}\nonumber \\
(\widehat{q}^{\mu}_{t};\widehat{p}^{\mu}_{t})&=&(q^{\mu,-q}_{t};p^{\mu,-p}_{t})
\textrm{ for } \tau_{2k+1}^{\mu}\leq t\leq \tau_{2k+2}^{\mu,-},
k=0,1,2,...\label{reflection_2nd}
\end{eqnarray}

Process $(\widehat{q}^{\mu}_{t};\widehat{p}^{\mu}_{t})$ is a process
with reflection and it can be seen that
$(\widehat{q}^{\mu}_{t};\widehat{p}^{\mu}_{t})$, which is the same
as
$(|q^{1,\mu}_{t}|,q^{2,\mu}_{t},\cdots,q^{r,\mu}_{t};\frac{d}{dt}|q^{1,\mu}_{t}|,\dot{q}^{2,\mu}_{t},\cdots,\dot{q}^{r,\mu}_{t})$,
and l.p.r.$(q^{\mu}_{t};p^{\mu}_{t})$ coincide.

In the figures below we give an illustration of the construction
of $(\widehat{q}^{1,\mu}_{t};\widehat{p}^{1,\mu}_{t})$. The first
figure illustrates with thick continuous and dotted lines
$\widehat{q}^{1,\mu}_{t}$ versus $t$. The continuous line is
$q^{1,\mu,q^{1}}_{t}$ versus $t$ and the dotted is
$q^{1,\mu,-q^{1}}_{t}$ versus $t$. The second figure illustrates
with thick continuous and dotted lines $\widehat{p}^{1,\mu}_{t}$
versus $t$. The continuous line is $p^{1,\mu,p^{1}}_{t}$ versus
$t$ and the dotted is $p^{1,\mu,-p^{1}}_{t}$ versus $t$.
\newpage
\begin{figure}[ht]
\begin{center}
\includegraphics[scale=0.4, width=5.5 cm, height=6 cm]{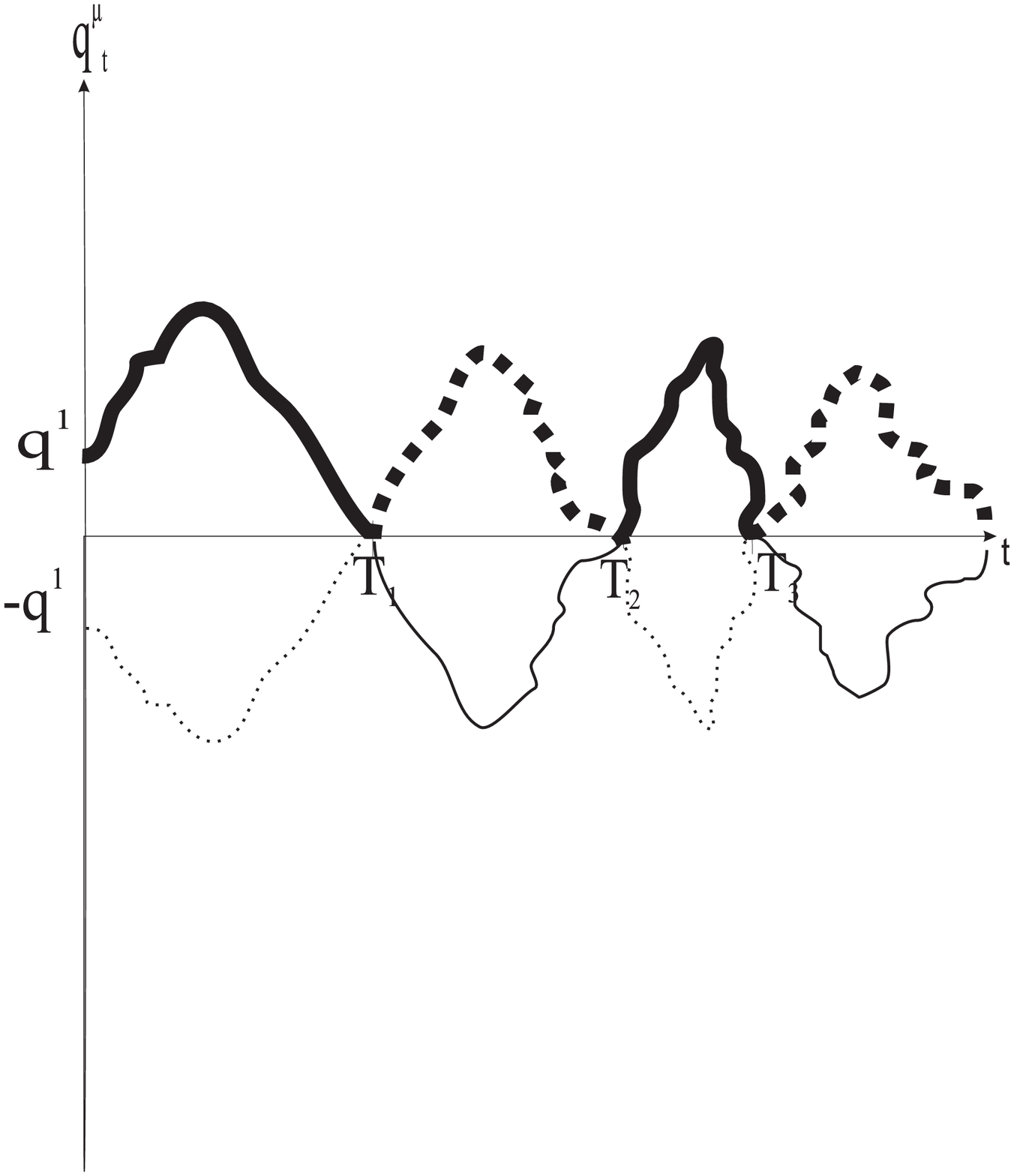}
\hspace{0.1cm}
\includegraphics[scale=0.4, width=5.5 cm, height=6 cm]{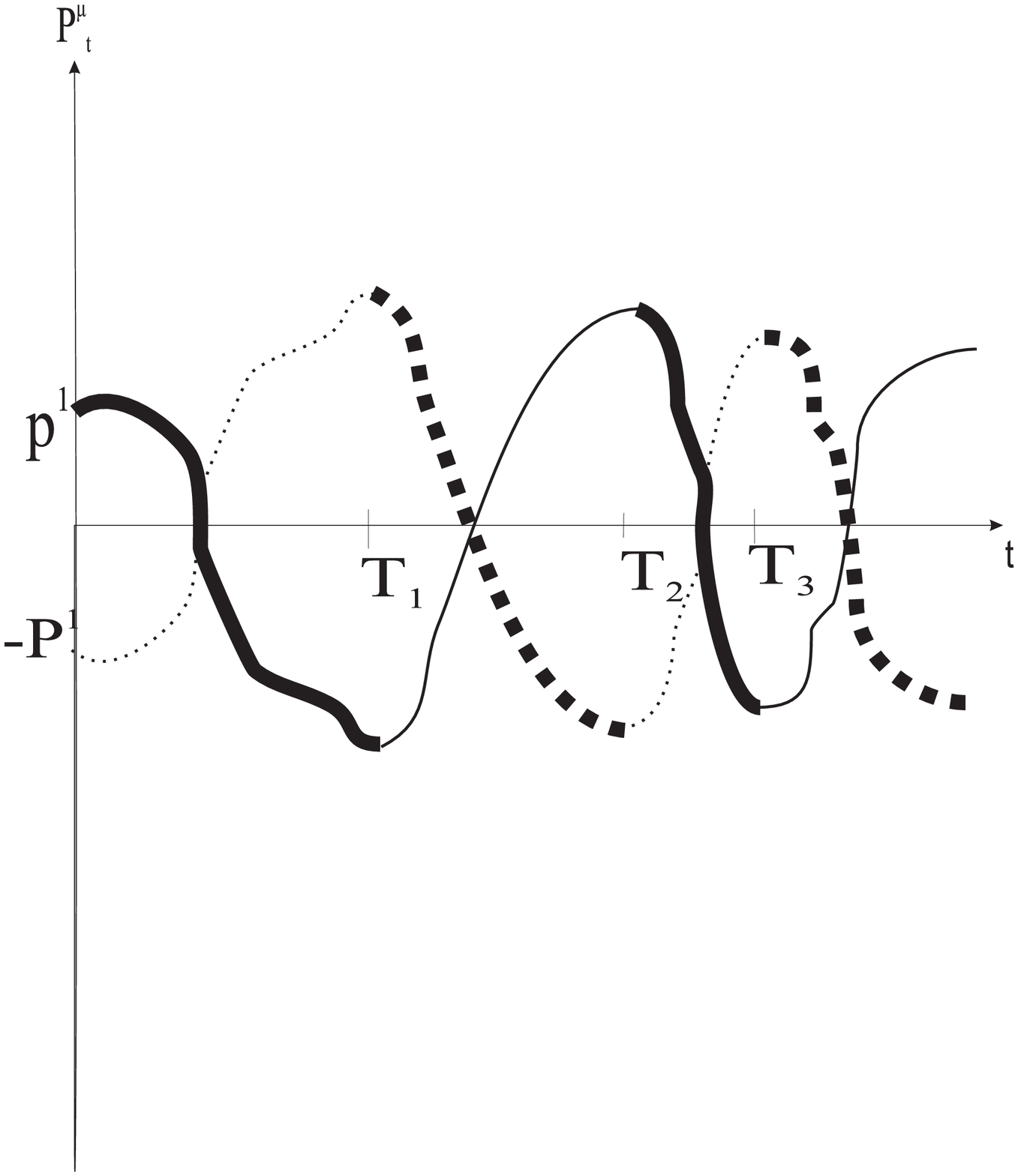}
\caption{Illustration of the process with reflection}
\end{center}
\end{figure}

\textbf{Lemma 2.3.} \textit{Let $T>0$. The Markov process
$(q^{\mu}_{t};p^{\mu}_{t})$ starting at a point $(q,p)$ different
from the origin $\textrm{O}=(0,...,0;0,...,0)$, that satisfies
system (\ref{sk1_sgn_withb}), does not reach the origin
$\textrm{O}$ in finite time T, i.e. $$P(\exists t\leq T
\hspace*{0.1cm}s.t.\hspace*{0.1cm}
(q^{\mu}_{t};p^{\mu}_{t})=\textrm{O})=0.$$}

\textbf{Proof.} We easily see that it is actually enough to
consider only $(q^{1,\mu}_{t};p^{1,\mu}_{t})$. Let $\delta\ll 1$
be a small number. Define the rectangle $\Delta=\{(q,p)\in
\mathbb{R}\times \mathbb{R}: |q|\leq \frac{\delta^{2}}{2},|p|\leq
\frac{\delta}{2}\}$ and suppose that the trajectory starts from
some point outside the rectangle $\Delta$, say from $(q,0)\in
\mathbb{R}^{2}\setminus \Delta$.  Let also $\chi_{\Delta}(x)$
denote the indicator function of the set $\Delta$. Then
$E^{(q,0)}\int_{0}^{T}\chi_{\Delta}(q_{s}^{1},p_{s}^{1})ds$ is the
expected value of the time ,during time $[0,T]$,  that the process
$(q_{t}^{1},p_{t}^{1})$ with initial point $(q,0)$ spends inside
the rectangle $\Delta$. If $b=0$ and $\Sigma$ is a matrix with
constant entries, $(q^{1}_{t},p^{1}_{t})$ is a Gaussian process.
One can write down its density explicitly (see equation
(\ref{sk1_sgn_withb})), which we denote by $\rho(\cdot)$, and
obtain the bound
\begin{equation}
E^{(q,0)}\int_{0}^{T}\chi_{\Delta}(q_{s}^{1},p_{s}^{1})ds=
\int_{\Delta}\int_{0}^{T}\rho(s,(q,0),y)ds dy\leq
A(T,q)\delta^{3}\label{upper_bound}
\end{equation}
where $A(T,q)$ is a constant that depends on $T$ and $q$. The
general case can be reduced to the case with $b=0$ and $\Sigma$
constant by an absolutely continuous change of measures in the
space of trajectories and by a random time change.

We will establish now a lower bound for the quantity
$E^{(q,0)}\int_{0}^{T}\chi_{\Delta}(q_{s}^{1},p_{s}^{1})ds$ under
the assumption that the process $(q_{t}^{1,\mu},p_{t}^{1,\mu})$
will reach $(0,0)$ before time $T$ with positive probability. This
will lead to a contradiction.

Again by Girsanov's theorem on the absolute continuity of measures
in the space of trajectories it is enough to consider the solution
of the following S.D.E:
\begin{eqnarray}
\dot{q}^{1}_{t}&=& p^{1}_{t}\nonumber \\
\dot{p}^{1}_{t}&=&\frac{1}{\mu}\sum_{j=1}^{r}\sigma_{1j}(q^{\mu}_{t})\dot{\overline{W}}_{t}^{j}\label{not_origin}\\
q^{1}_{0}&=&q^{1}, p^{1}_{0}=p^{1}, \nonumber
\end{eqnarray}
where
$\overline{W}_{t}^{j}=\int_{0}^{t}\textrm{sgn}(q_{u}^{1,\mu})dW_{u}^{j}$.

By the self similarity properties of the Wiener process one can
find a Wiener process $W_{t}^{1,*}$ such that
$\int_{0}^{t}\frac{1}{\mu}\sum_{j=1}^{r}\sigma_{1j}(q^{\mu}_{t})\dot{\overline{W}}_{t}^{j}=
W_{\theta(t)}^{1,*} $, where
$\theta(t)=\int_{0}^{t}\frac{1}{\mu^{2}}\alpha_{11}(q_{s}^{\mu})ds$
and $\alpha_{11}=\sum_{j,k=1}^{r}\sigma_{1j}\sigma_{1k}$. So
$\int_{0}^{t}\frac{1}{\mu}\sum_{j=1}^{r}\sigma_{1j}(q^{\mu}_{t})\dot{\overline{W}}_{t}^{j}$
can be obtained from  $W_{t}^{1,*}$ via a random time change.

By the law of iterated logarithm we get that for all $k\in[0,1]$
there exists a $t_{o}(k)$ small enough, such that
\begin{displaymath}
P(t^{\frac{1}{2}+k}\leq |W_{t}^{1,*}|\leq t^{\frac{1}{2}-k}
\textrm{ for } t\in [0,t_{o}(k)])\geq 1-k.
\end{displaymath}
Observe that if $t\in [0,t_{o}(k)]$ then $\theta(t) \in
[0,ct_{o}(k)]$, where $c=\frac{1}{\mu^{2}}\sup_{{x \in
\mathbb{R}}}|\alpha_{11}(x)|$. Define also
$t_{o}^{'}(k)=\min\{t_{o}(k),\frac{t_{o}(k)}{c}\}$. Then with
probability very close to 1, as $k \rightarrow 0$, and for all
$t\in [0,t_{o}^{'}(k)]$ it must hold that $|p_{t}^{1,\mu}|\leq
c_{1}t^{\frac{1}{2}-k}$ and $q_{t}^{1,\mu}=
\int_{0}^{t}p_{s}^{1,\mu} ds \leq \int_{0}^{t}
c_{1}s^{\frac{1}{2}-k}ds < 2c_{1}t^{\frac{3}{2}-k}$, for a
constant $c_{1}$.

Let $\tau$ be the first time, after the time that the Markov
process reached the origin, that it exits from the rectangle
$\Delta$, i.e. $\tau=\inf\{t>0:(q_{t}^{1},p_{t}^{1})\in
\mathbb{R}^{2}\setminus \Delta\}$. Then it follows that
\begin{equation}
E^{(q,0)}\int_{0}^{T}\chi_{\Delta}(q_{s}^{1},p_{s}^{1})ds>
E\{\tau\}\times P(\exists t\leq T
\hspace*{0.1cm}s.t.\hspace*{0.1cm} (q_{t}^{1,\mu};
p_{t}^{1,\mu})=(0,0))
\end{equation}

Define $\tau_{q}=\inf\{t>0:|q_{t}^{1,\mu}|>\frac{\delta^{2}}{2}\}$
and $\tau_{p}=\inf\{t>0:|p_{t}^{1,\mu}|>\frac{\delta}{2}\}$. By
the above bounds for $q_{t}^{1,\mu}$ and $p_{t}^{1,\mu}$ we get
that $\tau_{q}>c_{q}\delta^{\frac{4}{3}}$ and
$\tau_{p}>c_{p}\delta^{2}$, where $c_{q},c_{p}$ are some constants
independent of $\delta$. So the trajectory exits the rectangle
faster in the direction of $p$ than in the direction of $q$ and
the exit time is of order $\delta^{2}$. Therefore, by this and by
(\ref{upper_bound}), we have that
\begin{equation}
B\delta^{2}<
E^{(q,0)}\int_{0}^{T}\chi_{\Delta}(q_{s}^{1},p_{s}^{1})ds \leq
A\delta^{3},
\end{equation}
which cannot hold for constants A and B and small enough $\delta$.
So we have a contradiction and hence it is true that $P(\exists
t\leq T \hspace*{0.1cm}s.t.\hspace*{0.1cm}
(q_{t}^{1,\mu};p_{t}^{1\mu})=(0,0))=0$.
\begin{flushright}
$\Box$
\end{flushright}

\textbf{Theorem 2.4.}  \textit{We have the following two
statements:
\begin{enumerate}
\item {Let $T>0$. The Markov process
l.p.r.$(q^{\mu}_{t};p^{\mu}_{t})$ (with arbitrary $b$) does not
reach the origin $\textrm{O}=(0,...,0;0,...,0)$ in finite time
$T$, namely
$$P(\exists t\leq T \hspace*{0.1cm}s.t.\hspace*{0.1cm}
\textrm{l.p.r.}(q_{t}^{\mu};p_{t}^{\mu})=\textrm{O})=0.$$} \item
{The sequence of Markov times $\{\tau_{k}^{\mu}\}$ converges to
$+\infty$ as $k\rightarrow +\infty$, i.e. $$P(\lim_{k\rightarrow
+\infty}\tau_{k}^{\mu}=+\infty)=1$$}
\end{enumerate}}

\textbf{Proof.} The Langevin process with reflection,
l.p.r.$(q^{\mu}_{t};p^{\mu}_{t})$, coincides at any time $t$
either with $(q^{\mu,q}_{t};p^{\mu,p}_{t})$ or with
$(q^{\mu,-q}_{t};p^{\mu,-p}_{t})$. Therefore we have that:
\begin{eqnarray}
P(\exists t\leq T \hspace*{0.1cm}s.t.\hspace*{0.1cm}
\textrm{l.p.r.}(q_{t}^{\mu};p_{t}^{\mu})=\textrm{O})&\leq&
P(\exists t\leq T \hspace*{0.1cm}s.t.\hspace*{0.1cm}
(q^{\mu,q}_{t};p^{\mu,p}_{t})=\textrm{O})\nonumber\\
&+&P(\exists t\leq T \hspace*{0.1cm}s.t.\hspace*{0.1cm}
(q^{\mu,-q}_{t};p^{\mu,-p}_{t})=\textrm{O}).\nonumber
\end{eqnarray}

Hence Lemma 2.3 implies that
\begin{displaymath}
P(\exists t\leq T \hspace*{0.1cm}s.t.\hspace*{0.1cm}
\textrm{l.p.r.}(q_{t}^{\mu};p_{t}^{\mu})=\textrm{O})=0.
\end{displaymath}

Part (ii) is an easy consequence of part (i). It is easy to see
that $\{\tau_{k}^{\mu}\}$ is an unbounded, strictly increasing
sequence of Markov times. Indeed, if on the contrary we assume
that there exists a $N$ such that $\tau_{k}^{\mu}\leq N$ for all
$k$ with positive probability, then the trajectories of
l.p.r.$(q^{\mu}_{t};p^{\mu}_{t})$ will have limit points. The only
possible limit point however is the origin $(0,...,0;0,...,0)$.
But by part (i) the probability that within any time $T$ the
trajectory will reach the origin is 0. So $\{\tau_{k}^{\mu}\}$ is
an unbounded strictly increasing sequence of Markov times.
Therefore we have that $P(\lim_{k\rightarrow
+\infty}\tau_{k}^{\mu}=+\infty)=1$.

\begin{flushright}
$\Box$
\end{flushright}

Therefore the Langevin process with reflection has only finitely
many jumps in any time interval $[0,T]$ with probability 1. Hence
our definition for the Langevin process with reflection is
correct.

\section{The Skorohod reflection problem}

The convergence of the Langevin process with reflection that will
be presented in section 4 relies on results about solutions of the
Skorohod reflection problem, proven in \cite{C2} and \cite{T1}.

Let us first recall that $D=\mathbb{R_{+}}\times
\mathbb{R}^{r-1}$, $\partial D=\partial \mathbb{R_{+}}\times
\mathbb{R}^{r-1}$ and let $N(q)$ be the set of inward normals at
$q\in\partial D$. Denote also by $\mathbb{D}(\mathbb{R}_{+},D)$
the space of \textit{cadl\'{a}g} (right continuous with left
limits) functions with values in $D$, endowed with the Skorohod
topology and by $\mathbb{B.V.}(\mathbb{R}_{+},D)$ the set of
\textit{cadl\'{a}g} functions with bounded variation and values in
$D$.

\vspace{0.2cm}

\textbf{Definition 3.1.}\textit{ Let $w$ be a function in
$\mathbb{D}(\mathbb{R}_{+},\mathbb{R}^{r})$ such that $w(0)\in D$.
We say that the pair $(q,\phi)$ with $q\in
\mathbb{D}(\mathbb{R}_{+},D)$, $\phi \in
\mathbb{B.V.}(\mathbb{R}_{+},\mathbb{R}^{r})$ is a solution to the
Skorohod problem for $(D,N,w)$ if}
$$q_{t}=w_{t}+\phi_{t}$$
$$\phi_{t}=\int_{0}^{t}\nu(s)d|\phi|_{s}, \nu(s)\in N(q_{s}), d|\phi|-a.e.$$
$$d|\phi|({t:q_{t}\in D})=0,$$

\textit{where $|\phi|$ denotes the total variation of $\phi$ and
is called the local time of the solution.}

\vspace{0.2cm}

The following theorem characterizes the continuity properties of
solutions of the Skorohod reflection problem.

\vspace{0.2cm}

\textbf{Theorem 3.2.} \textit{Let $W$ be a compact subset of
$\mathbb{D}(\mathbb{R}_{+},\mathbb{R}^{r})$ in the Skorohod
topology such that $w(0)\in D$ for every $w \in W$. Moreover let
$\mathcal{Q}$ be the set of $(q,\phi,|\phi|,w)\in
\mathbb{D}(\mathbb{R}_{+},D)\times
\mathbb{B.V.}(\mathbb{R}_{+},\mathbb{R}^{r}) \times
\mathbb{B.V.}(\mathbb{R}_{+},\mathbb{R}_{+}) \times
 \mathbb{D}(\mathbb{R}_{+},\mathbb{R}^{r})$ such that $(q,\phi)$
 is the solution to the Skorohod problem for $(D,N,w)$ for some $w \in
 W$ and $q$ is continuous. The set $D$ is convex and so $\mathcal{Q}$ is
 a relatively compact subset of
 $\mathbb{D}(\mathbb{R}_{+},\mathbb{R}^{3r+1})$ in the Skorohod
 topology and for every accumulation point of $(q,\phi,|\phi|,w)$
 in $\mathcal{Q}$ we have that $(q,\phi)$ is a solution to the
 Skorohod problem for $(D,N,w)$.}

\textbf{Proof.} This is a special case of theorem 3.2 in
\cite{C1}.
\begin{flushright}
$\Box$
\end{flushright}

\section{Convergence of the Langevin process with reflection}

In this section we consider the limit of l.p.r.$(q_{t}^{\mu})$ as
$\mu\rightarrow 0$ when the diffusion matrix is the unit matrix.
Below we will assume that $t\leq T$, where $T$ ia s positive real
number.

Consider the stochastic process $(q_{t}^{\mu};p_{t}^{\mu})$ in
$D\times \mathbb{R}^r$, which satisfies the following system of
S.D.E.'s:

\begin{eqnarray}
\dot{q}^{\mu}_{t}&=& p^{\mu}_{t}\nonumber \\
\mu\dot{p}^{\mu}_{t}&=&-p^{\mu}_{t}+b(q^{\mu}_{t})+ \dot{W}_{t}+\nu(q^{\mu}_{t})\cdot\dot{\Psi}^{\mu}_{t}\label{langevin_jump1}\\
q^{\mu}_{0}&=&q_{0}, p^{\mu}_{0}=p_{0}, \nonumber
\end{eqnarray}

where $q_{t}^{\mu}=(q_{t}^{1,\mu},\cdots,q_{t}^{r,\mu})^{'}$,
$p_{t}^{\mu}=(p_{t}^{1,\mu},\cdots,p_{t}^{r,\mu})^{'}$,
$W_{t}=(W_{t}^{1},\cdots,W_{t}^{r})^{'}$, $\nu(q)$ denotes the
unit inward normal to $D$ at $q\in\partial D$,
$b(q)=(b_{1}(q),...,b_{r}(q))^{'}$ and
$\Psi^{\mu}_{t}=\mu\sum_{s\leq
t}(-2p_{s-}^{\mu}\cdot\nu(q_{s}^{\mu}))\cdot\chi_{\partial
D}(q_{s}^{\mu})$. It is easy to see that (\ref{langevin_jump1}) is
pathwise equivalent to the Langevin process with reflection in
$D\times \mathbb{R}^r$ of Definition 2.1. and so it admits a
unique weak solution.

We will follow the method introduced in \cite{C1}. The main idea
is to represent $q^{\mu}$ as the first component of a solution to
the Skorohod problem for $(D,N, H^{\mu}+ X^{\mu})$, where $H^{\mu}
+ X^{\mu}$ is a semimartingale. The family $\{H^{\mu}+ X^{\mu}\}$
turns out to be tight and this enables us to use Theorem 3.2 to
conclude that the family $\{q^{\mu}\}$ is tight as well.

We can suppose that there is a unique underlying complete
probability space $(\Omega, \mathbb{F},P)$. Let
$\mathbb{\widehat{F}}$ denote the the $\sigma-$algebra of
$\mathbb{F}$ of sets with $P-$ measure $0$ or $1$ and define the
filtration
$$\mathbb{F}^{\mu}_{t}=\mathbb{\widehat{F}}\cup \sigma((q_{s}^{\mu};p_{s}^{\mu}), s\leq t).$$

\textbf{Lemma 4.1.} \textit{For every $\mu$ the pair of stochastic
processes $(q^{\mu}_{\cdot},\Phi^{\mu}_{\cdot})$, where
\begin{equation}
\Phi^{\mu}_{t}=\int_{0}^{t}\nu(q_{s}^{\mu})d\Psi^{\mu}_{t}
\label{local_time1}
\end{equation}
is an almost surely solution to the Skorohod reflection problem
for $(D,N,H^{\mu}+ X^{\mu})$, where
\begin{eqnarray}
H^{\mu}_{t}&=&q_{0}+\mu p_{0}- \mu p_{t}^{\mu}\nonumber \\
X^{\mu}_{t}&=&\int_{0}^{t}b(q_{s}^{\mu})ds
+W_{t}\label{semimartingale1}
\end{eqnarray}}

\textbf{Proof.} Consider the integral form of
(\ref{langevin_jump1}). Taking into account that
$\int_{0}^{t}p_{s}^{\mu}ds=q_{t}^{\mu}-q_{0}$ and solving for
$q_{t}^{\mu}$ we see that:
\begin{displaymath}
q_{t}^{\mu}=H^{\mu}_{t}+ X^{\mu}_{t}+ \Phi^{\mu}_{t}
\end{displaymath}
Then $(q^{\mu},\Phi^{\mu})$ verifies Definition 3.1 with
probability 1.
\begin{flushright}
$\Box$
\end{flushright}

\textbf{Lemma 4.2.} \textit{For every $T>0$ we have that
$\lim_{\mu\rightarrow 0}E[\sup_{t\leq T}|\mu p_{t}^{\mu}|^{2}]=0$.
} \vspace{0.1cm}

\textbf{Proof.} Assume first that $b=0$.
 Consider equations (\ref{langevin_jump1}) and apply the
 It$\widehat{o}$
 formula for semimartingales to the function $f(q,p)=|p|^{2}$ for
 every pair of times $s,t$ such that $0\leq s\leq t \leq
 T$. Doing that we get

 \begin{eqnarray}
|p^{\mu}_{t}|^{2}&=&|p^{\mu}_{s}|^{2}
-\frac{2}{\mu}\int_{s}^{t}|p^{\mu}_{u}|^{2}du
+\frac{2}{\mu}\int_{s}^{t}p^{\mu}_{u} \cdot dW_{u}
+\frac{1}{\mu^{2}}r(t-s) \label{p_square_ito}
\end{eqnarray}

It is interesting to observe that the local time $\Psi_{t}^{\mu}$
does not appear above. This comes from the fact that under elastic
reflection  $|p^{\mu}_{t}|^{2}= |p^{\mu}_{t-}|^{2}$ for every
$t>0$.

Consider now a constant $c>0$ and functions $x,g\in
\mathbb{D}([0,T],\mathbb{R})$ with $g(0)=0$ such that:

\begin{equation}
x_{t}\leq x_{s}- c\int_{s}^{t}x_{u}du+g_{t}-g_{s},\hspace{0.1cm}
0\leq s\leq t\leq T \label{inequality1}
\end{equation}

Then one can easily see that
\begin{equation}
x_{t}\leq
e^{-ct}(x_{0}+g_{t})+c\int_{0}^{t}e^{-c(t-u)}(g_{t}-g_{u})du,\hspace{0.1cm}
0\leq t\leq T \label{inequality2}
\end{equation}

By taking expected value to (\ref{p_square_ito}) and applying
(\ref{inequality2}) with $c=\frac{2}{\mu}$,
$g_{t}=\frac{1}{\mu^{2}}rt$ and $x_{t}=|p^{\mu}_{t}|^{2}$, we get

\begin{eqnarray}
E|p^{\mu}_{t}|^{2}&\leq&
e^{-\frac{2}{\mu}t}(|p|^{2}+\frac{1}{\mu^{2}}rt)+
\frac{2}{\mu^{3}}\int_{0}^{t}e^{-\frac{2}{\mu}(t-u)}r(t-u)du\nonumber\\
&=&e^{-\frac{2}{\mu}t}|p|^{2}+ \frac{r}{\mu^{2}}
(\frac{\mu}{2}-\frac{\mu}{2}e^{-\frac{2t}{\mu}})\label{Ep_square_ito}
\end{eqnarray}

This implies the statement of the Lemma for $b=0$. The general
case can be reduced to the case with $b=0$ by an absolutely
continuous change of measures in the space of trajectories.
\begin{flushright}
$\Box$
\end{flushright}

The following two theorems are restatements of theorems 3.8.6 and
3.10.2 respectively of \cite{EK1}.

\vspace{0.2cm}

\textbf{Theorem 4.3.} \textit{Let $\{Y^{n}\}$ be a family of processes
with sample paths in $\mathbb{D}(\mathbb{R}_{+},D)$. Assuming that
for every $\epsilon>0$ and rational $t\geq 0$ there exist a
compact set $\Gamma(\epsilon,t)\subset D$ such that $\liminf_{n}P
(Y^{n}(t)\in \Gamma(\epsilon,t))\geq 1-\epsilon$, then the
following are equivalent
\begin{enumerate}
\item{$\{Y^{n}\}$ is relatively compact.} \item{For each $T>0$,
there exists $\beta>0$ and a family of nonnegative random
variables $\{\gamma^{n}(\delta),0<\delta<1\}$ satisfying
      $$E(|Y^{n}(t+u)-Y^{n}(t)|^{\beta} \hspace{0.1cm}|\mathbb{F}^{n}_{t})\leq E(\gamma^{n}(\delta)|\mathbb{F}^{n}_{t}),$$
      for $t\in [0,T]$ and $u\in [0,\delta]$ and in addition
      $\lim_{\delta\rightarrow
      0}\limsup_{n}E(\gamma^{n}(\delta))=0
      $.}
\end{enumerate}}

\vspace{0.2cm}

\textbf{Theorem 4.4.} \textit{Let $\{Y^{n}\}$ and $Y$ be processes
with sample paths in $\mathbb{D}(\mathbb{R}_{+},D)$ such that
$Y_{n}$ converges in distribution to $Y$. Then $Y$ is almost
surely continuous if and only if $\int_{0}^{\infty}
e^{-u}[\sup_{0\leq t \leq u}|Y^{n}(t)-Y^{n}(t-)|\wedge
1]du\Rightarrow 0$. } \vspace{0.2cm}

The following lemma shows that the family $\{H^{\mu}+X^{\mu}\}$ is
tight in the Skorohod topology.

\vspace{0.2cm}

 \textbf{Lemma 4.5.} \textit{The family $\{H^{\mu}+X^{\mu}\}$ defined in (\ref{semimartingale1}) is relatively
compact and all of its accumulation points are continuous.}

\textbf{Proof.} It is easily seen that $\{X^{\mu}\}$ is relatively
compact and that all of its accumulation points are continuous.

Now Lemma 4.2 suggests that:
\begin{eqnarray}
\lim_{\mu\rightarrow 0}E[\sup_{t\leq T}|H_{t}^{\mu}|^{2}]\leq
c\label{rel_compactness_of_H_1}\\
\lim_{\mu\rightarrow 0}E[\sup_{|t-s|\leq
\delta}|H_{t}^{\mu}-H_{s}^{\mu}|]\leq
c_{1}\delta,\label{rel_compactness_of_H_2}
\end{eqnarray}
where $c,c_{1}$ are positive constants independent of $\mu$.

Chebychev's inequality and (\ref{rel_compactness_of_H_1}) imply
that $$\liminf_{n\rightarrow \infty}P (|H^{1/n}(t)|\leq
\lambda)\geq 1-\frac{c}{\lambda^{2}}.$$ Therefore by this and
(\ref{rel_compactness_of_H_2}), Theorem 4.3. gives us that
$\{H^{\mu}\}$ is relatively compact. Lastly
(\ref{rel_compactness_of_H_2}) and Theorem 4.4 implies that all
its accumulation points are continuous.

\begin{flushright}
$\Box$
\end{flushright}

\textbf{Theorem 4.6.} \textit{The family
$\{(q^{\mu},\Phi^{\mu},\Psi^{\mu},H^{\mu},X^{\mu})\}$ is
relatively compact in
$\mathbb{D}(\mathbb{R}_{+},\mathbb{R}^{4r+1})$. }

\textbf{Proof.} It follows from Lemma 4.5 and Theorem 3.2.

\begin{flushright}
$\Box$
\end{flushright}

Now that tightness has been established we will proceed with the
identification of the stochastic differential equation with
reflection that describes the behavior of $q^{\mu}$ as
$\mu\rightarrow 0$.

Consider the following S.D.E. with reflection:
\begin{equation}
q_{t}=q_{0}+\int_{0}^{t}b(q_{s})ds +W_{t}+\Phi_{t}\label{limiting_process1}\\
\end{equation}
where $\Phi_{t}=\int_{0}^{t}\nu(q_{s})d|\Phi|_{s},\nu(s)\in
N(q_{s})$ and $d|\Phi|(\{t:q_{t}\in D\})=0$. It is known that
(\ref{limiting_process1}) has a unique weak solution $(q,\Phi)$
(\cite{AO1}).

\vspace{0.2cm}

\textbf{Theorem 4.7.} \textit{The family
$\{(q^{\mu},\Phi^{\mu})\}$ converges in distribution to the unique
solution $(q,\Phi)$ of (\ref{limiting_process1}).}

\textbf{Proof.} By Theorem 4.6. we have that the five-tuple
$\{(q^{\mu},\Phi^{\mu},H^{\mu},X^{\mu},W)\}$ is relatively compact
in $\mathbb{D}(\mathbb{R}_{+},\mathbb{R}^{5r})$. Hence it (or a
subsequence) converges in distribution to a stochastic process
$\{(q,\Phi,H,X,W)\}$. By the Skorohod representation theorem, one
can find a probability space
$(\widetilde{\Omega},\widetilde{\mathbb{F}},\widetilde{P})$ and
realizations
$\{(\widetilde{q}^{\mu},\widetilde{\Phi}^{\mu},\widetilde{H}^{\mu},\widetilde{W}^{\mu})\}$
and
$\{(\widetilde{q},\widetilde{\Phi},\widetilde{H},\widetilde{X},\widetilde{W})\}$
of $\{(q^{\mu},\Phi^{\mu},H^{\mu},X^{\mu},W)\}$ and
$\{(q,\Phi,H,X,W)\}$ respectively such that
$\{(\widetilde{q}^{\mu},\widetilde{\Phi}^{\mu},\widetilde{H}^{\mu},\widetilde{X}^{\mu},\widetilde{W}^{\mu})\}$
converges $\widetilde{P}$-almost surely to
$\{(\widetilde{q},\widetilde{\Phi},\widetilde{H},\widetilde{X},\widetilde{W})\}$.
Therefore by Theorem 3.2. $(\widetilde{q}, \widetilde{\Phi})$ is a
solution to the Skorohod problem for
$(D,N,\widetilde{H}+\widetilde{X})$ $\widetilde{P}-$almost surely.

Now by the convergence of $\widetilde{q}^{\mu}$ to $\widetilde{q}$
we get that $\widetilde{X}$ must be given by:
\begin{displaymath}
\widetilde{X}_{t}=\int_{0}^{t}b(\widetilde{q}_{s})ds +
\widetilde{W}_{t}
\end{displaymath}

Finally Lemma 4.2 and its proof imply that
$\widetilde{H}_{t}=q_{0}$.
\begin{flushright}
$\Box$
\end{flushright}

We would like to note here that one could prove the convergence in
distribution of the Langevin procces with reflection to the
corresponding diffusion process with reflection using the
Smoluchowski-Kramers approximation. However the beauty and
generality of the results of \cite{C2} resulted in using the
method that was presented here.

\section*{Acknowledgments}
I would like to thank my advisor Mark Freidlin for posing the
problem and for his valuable help and Dwijavanti Athreya, Hyejin
Kim and James (J.T.) Halbert for their helpful suggestions.


\end{document}